\documentclass{article}

\usepackage{amsmath}
\usepackage{amssymb}
\usepackage{amsfonts}
\usepackage{latexsym}
\usepackage{enumerate}

\newtheorem{llemma}{Lemma}[section]

\newtheorem{exmp}[llemma]{Example}
\newtheorem{ttheorem}[llemma]{Theorem}

\newtheorem{defn}[llemma]{Definition}
\newtheorem{rem}[llemma]{Remark}

\newtheorem{key}[llemma]{Keyword}

\begin{document}

\title{\textbf{Convex and Concave Soft Sets and Some Properties}}
\author{Irfan Deli \\ \\
             Department of Mathematics, Faculty of Arts and
          Sciences,\\
 Kilis 7 Aral{\i}k University, 79000 Kilis, Turkey, \\
              irfandeli@kilis.edu.tr
       }

\maketitle


\begin{abstract}
In this study, after given the definition of soft sets and their
basic operations we define convex soft sets which is an important
concept for operation research, optimization and related problems.
Then, we define concave soft sets and give some properties for the
concave sets. For these, we will use definition and properties of
convex-concave fuzzy sets in literature. We also give different some
properties for the convex and concave soft sets.
\begin{key}
Soft sets, fuzzy set,  $ \alpha $-inclusion, Convex soft sets,
concave soft sets.
\end{key}
\end{abstract}


\section{Introduction}
Many fields deal with the uncertain data which may not be
successfully modeled by the classical mathematics, probability
theory, fuzzy sets \cite{zad-65}, rough sets \cite{paw-82}, and
other mathematical tools. In 1999, Molodtsov \cite{mol-99} proposed
a completely new approach so-called {\em soft set theory} that is
more universal for modeling vagueness and uncertainty.

After definitions of the operations of soft sets
\cite{ali-09,cag-09a,maj-03,aslý-11a}, the properties and
applications on the soft set theory have been studied increasingly
(e.g. \cite{ali-09,cag-11d,cag-09b,deli-13a,deli-13b}). The
algebraic structure of soft set theory has also been studied
increasingly (e.g. \cite{acar-10,akt-07,jun-08b,zhan-10}). In recent
years, many interesting applications of soft set theory have been
expanded by embedding the ideas of fuzzy sets (e.g. \cite{ ali-11,
deli-12a,deli-12b,fen-08b,zhan-10,roy-07,maj-01a}), rough sets (e.g.
\cite{ali-11, fen-11, fen-10}) and intuitionistic fuzzy sets (e.g.
\cite{jia-11,maj-01c, maj-01b,muk-08}).

Convex (concave) and convex (concave )fuzzy sets play important
roles in optimization theory. Various definitions of convex fuzzy
and concave fuzzy sets have appeared in the literature. A
significant definition of convex fuzzy sets introduced by Zadeh
\cite{zad-65} and concave fuzzy sets introduced by Chaudhuri
 \cite{cha-92}. Then, concavoconvex fuzzy sets proposed by Sarkar
\cite{sar-96}, which is convex and concave fuzzy sets together
conceived by combining. Moreover, some properties of the
concavoconvex fuzzy set are given by the author.  Subsequently,
works on convex (concave )fuzzy sets in theories and applications
has been progressing rapidly by \cite{amm-99,wan-09,yan-95}.

In this work, we introduce the soft version of the Zadeh's
definition of fuzzy convex set, and we will call this generalization
an soft convex sets. Also, we will introduce the definition of soft
concave set and we study some desired properties.

The plan of the paper is as follows. In Section 2, we give some
notations, definitions used throughout the paper. In Section 3, we
define soft convex sets and soft concave sets and then we show some
properties.
\section{Preliminary}\label{ss}

In this section, we present the basic definitions and results of
soft set theory \cite{cag-09a}. More detailed explanations related
to this subsection may be found in \cite{cag-09a, maj-03, mol-99}.
\begin{defn}\cite{mol-99}
Let $U$ be a universe, $P(U)$ be the power set of $U$ and $E$ be a
set of parameters that are describe the elements of $U$. A soft set
$S$ over $U$ is a set defined by a set valued function $S$
representing a mapping
$$
f_S: E\to P(U)
$$
It is noting that the soft set is a parametrized family of subsets
of the set $U$, and therefore it can be written a set of ordered
pairs
$$
S= \{(x, f_S(x)): x\in E\}
$$
\end{defn}
Here, $f_S$ is called approximate function of the soft set $S$ and
$f_S(x)$ is called $x$-approximate value of $x \in E$. The subscript
$S$ in the $f_S$ indicates that $f_S$ is the approximate function of
$S$.

Generally, $f_S$, $f_T$, $f_V$, ... will be used as an approximate
functions of $S$, $T$, $V$, ..., respectively.

Note that if $f_S(x)=\emptyset$, then the element $(x, f_S(x))$ is
not appeared in $S$.
\begin{exmp}\label{ex-soft}
Suppose that $U=\{u_1,u_2,u_3,u_4,u_5,u_6\}$ is the universe and
$E=\{x_1,x_2,x_3,x_4,x_5\}$ is the set of parameters. If $f_S(x_1)=
\{u_1, u_2,u_3\}$, $f_S(x_2)= \{u_1,u_4 \}$, $ f_S(x_3)= \emptyset$,
$ f_S(x_4)= U$ then the soft-set $S$ is written by
$$
S=\{(x_1, \{u_1, u_2,u_3\}),(x_2, \{u_1,u_4\}), (x_4,U)\}
$$

By using same parameter set $E$, we can construct a soft set $T$
different then $S$. Assume that $f_T(x_1)= \{u_1, u_2,u_6\}$,
$f_T(x_2)= \{u_1, u_2,u_3 \}$, $ f_T(x_3)=\{u_1, u_2\}$, $ f_T(x_4)=
\emptyset$, $ f_T(x_5)= \{u_6\}$ then the soft-set $T$ is written by
$$
T=\{(x_1, \{u_1, u_2,u_6\}), (x_2, \{u_1, u_2,u_3\}), (x_3, \{u_1,
u_2\}), (x_5,\{u_6\}) \}
$$
\end{exmp}
\begin{defn}\cite{cag-09a} Let $S$ and $T$ be two soft sets. Then,
\begin{enumerate}[a)]
    \item
If $f_S(x)=\emptyset$ for all $x\in E$, then $S$ is called a empty
soft set, denoted by $S_\Phi$.
    \item
If $f_S(x)\subseteq f_T(x)$ for all $x\in E$, then $S$ is a soft
subset of $T$, denoted by $S \tilde{\subseteq} T$.
\end{enumerate}
\end{defn}
\begin{defn}\cite{cag-09a} Let $S$ and $T$ be two soft sets. Then,
\begin{enumerate}[a)]
\item
Complement of $S$ is denoted by $S^{\tilde c}$. Its approximate
function $f_{S^{\tilde c}}$ is defined by
$$
f_{S^{\tilde c}}(x)= U\setminus f_S(x)\quad\textrm{for all } x \in E
$$
\item
Union of $ S$ and $T$ is denoted by $S \tilde{\cup} T$. Its
approximate function $f_{S \widetilde{\cup} T}$ is defined by
$$
f_{S \widetilde{\cup} T}(x)= f_S(x)\cup f_T(x)\quad\textrm{for all }
x \in E.
$$

\item Intersection of $S$ and $T$ is
denoted by $S\tilde{\cap}T$. Its approximate function $f_{S
\widetilde{\cap} T}$ is defined by
$$
f_{S \widetilde{\cap} T}(x)= f_S(x)\cap f_T(x)\quad\textrm{for all }
x \in E.
$$
\end{enumerate}
\end{defn}
\begin{defn}\cite{cag-12e}  Let S be a soft set over U and $ \alpha $ be
a subset of U. Then, $ \alpha $-inclusion of the soft set S, denoted
by $S^\alpha$, is defined as

$$S^\alpha =\{ x\in E : f_S(x)\supseteq \alpha \}$$
\end{defn}

\section{Convex and concave soft sets}\label{ss}
In this section,  we define convex soft sets and concave soft sets
and then give desired some properties. Some of it is quoted from
example in \cite{amm-99,wan-09,sar-96,yan-95,zad-65}.

In this paper, E will denote the $n$-dimensional Euclidean space
$R^n$ and U denotes the arbitrary set.
\begin{defn}
The soft set S on E is called a convex soft set if

$$f_S(ax + (1 - a)y) \supseteq f_S(x)\cap f_S(y)$$ for every $x,y\in E$ and $a\in I$.
\end{defn}
%
\begin{ttheorem}
$S\cap T$ is a convex soft set when both S and T are convex soft
sets.
\end{ttheorem}
\emph{\textbf{Proof:}} Suppose that there exist $x, y\in E$ and
$a\in I$ and  $W= S \cap T$. Then,
\begin{equation}\label{6}
 f_W(ax +
(1 - a)y) =  f_S(ax + (1 - a)y) \cap  f_T(ax + (1 - a)y)
\end{equation}
Now, since S and T convex,
\begin{equation}\label{7}
 f_S(ax +
(1 - a)y) \supseteq f_S(x ) \cap  f_S(y)
\end{equation}
\begin{equation}\label{8}
 f_T(ax +
(1 - a)y) \supseteq f_T(x) \cap  f_T(y)
\end{equation}
and hence,
\begin{equation}\label{9}
 f_W(ax + (1 - a)y)\supseteq (f_S(x ) \cap  f_S(y))\cap (f_T(x) \cap
 f_T(y))
\end{equation}
and thus
\begin{equation}\label{10}
 f_W(ax +
(1 - a)y)\supseteq  f_W(x) \cap  f_W(y)
\end{equation}
\begin{rem}
If $\{S_i: i\in \{1,2,...\}\}$ is any family of convex soft sets,
then the intersection $\cap_{i\in I}S_i$ is a convex soft set.
\end{rem}
\begin{ttheorem}
The union of any family $\{S_i: i\in I=\{1,2,...\}\}$ of convex soft
sets is not necessarily a convex soft set.
\end{ttheorem}
\emph{\textbf{Proof:}}The proof is straightforward.

\begin{ttheorem} S is a convex soft set on E iff for every $\beta \in [0,1]$ and $\alpha \in P(U)$, $S^\alpha$ is a
convex set on E.
\end{ttheorem}
\emph{\textbf{Proof:}}$\Rightarrow$ Assume that S is a convex soft
set. If $x_1, x_2\in E$ and $\alpha \in P(U)$, then
$f_S(x_1)\supseteq \alpha$ and $f_S(x_2)\supseteq \alpha$. It
follows from the convexity of S that
$$
f_S(\beta x_1 +(1-\beta)x_2)\supseteq f_S(x_1)\cap f_S(x_2)
$$
and thus $S^\alpha$ is a convex set.

$\Leftarrow$ Assume that $S^\alpha$ is a convex set for every
$\beta\in [0, 1]$. Especially, for $x_1, x_2\in E,$ $S^\alpha$ is
convex for $\alpha =f_S(x_1)\cap f_S(x_2)$.

Since $f_S(x_1) \supseteq\alpha$ and $f_S(x_2) \supseteq\alpha$, we
have $x_1 \in S^\alpha$ and $x_2 \in S^\alpha$, whence $\beta x_1
+(1-\beta)x_2\in S^\alpha$. Therefore, $f_S(\beta x_1
+(1-\beta)x_2)\supseteq\alpha=  f_S(x_1)\cap f_S(x_2)$, which
indicates S is a convex soft set on X.
\begin{defn}
The soft set S on E is called a concave soft set if
$$f_S(ax + (1 - a)y) \subseteq f_S(x)\cup f_S(y)$$ for every $x,y\in E$ and $a\in I$.
\end{defn}
\begin{ttheorem}
$S\cup T$ is a concave soft set when both S and T are concave soft
sets.
\end{ttheorem}
\emph{\textbf{Proof:}} Suppose that there exist $x, y\in E$ and
$a\in I$ and  $W= S \cup T$. Then,
\begin{equation}\label{z6}
 f_W(ax +
(1 - a)y) =  f_S(ax + (1 - a)y) \cup  f_T(ax + (1 - a)y)
\end{equation}
Now, since S and T concave,
\begin{equation}\label{z7}
 f_S(ax +
(1 - a)y) \subseteq f_S(x ) \cup  f_S(y)
\end{equation}
\begin{equation}\label{z8}
 f_T(ax +
(1 - a)y) \subseteq f_T(x) \cup  f_T(y)
\end{equation}
and hence,
\begin{equation}\label{z9}
 f_W(ax + (1 - a)y)\subseteq (f_S(x ) \cup  f_S(y))\cup (f_T(x) \cup
 f_T(y))
\end{equation}
and thus
\begin{equation}\label{z10}
 f_W(ax +
(1 - a)y)\subseteq  f_W(x) \cup  f_W(y)
\end{equation}
\begin{rem}
If $\{S_i: i\in \{1,2,...\}\}$ is any family of concave soft sets,
then the union $\cup_{i\in I}S_i$ is a concave soft set.
\end{rem}
\begin{ttheorem}
$S\cap T$ is a concave soft set when both S and T are concave soft
sets.
\end{ttheorem}
\emph{\textbf{Proof:}} Suppose that there exist $x, y\in E$ and
$a\in I$ and  $W= S \cap T$. Then,
\begin{equation}\label{t6}
 f_W(ax +
(1 - a)y) =  f_S(ax + (1 - a)y) \cap  f_T(ax + (1 - a)y)
\end{equation}
Now, since S and T concave,
\begin{equation}\label{t7}
 f_S(ax +
(1 - a)y) \subseteq f_S(x ) \cup  f_S(y)
\end{equation}
\begin{equation}\label{t8}
 f_T(ax +
(1 - a)y) \subseteq f_T(x) \cup  f_T(y)
\end{equation}
and hence,
\begin{equation}\label{t9}
 f_W(ax + (1 - a)y)\subseteq (f_S(x ) \cup  f_S(y))\cap (f_T(x) \cup
 f_T(y))
\end{equation}
and thus
\begin{equation}\label{t10}
 f_W(ax +
(1 - a)y)\subseteq  f_W(x) \cap  f_W(y)
\end{equation}
\begin{rem}
The intersection of any family $\{S_i: i\in I=\{1,2,...\}\}$ of
concave soft sets is concave soft set.
\end{rem}
\begin{ttheorem}
$S^c$ is a concave soft set when S is a convex soft sets.
\end{ttheorem}
\emph{\textbf{Proof:}} Suppose that there exist $x, y\in E$, $a\in
I$ and S be a convex soft set.

Then, since S is convex,
\begin{equation}\label{x0}
f_S(ax + (1 - a)y) \supseteq f_S(x)\cap f_S(y)
\end{equation}
or
\begin{equation}\label{x1}
U\setminus f_S(ax + (1 - a)y) \subseteq U\setminus \{f_S(x)\cap
f_S(y)\}
\end{equation}
If $f_S(x)\supset f_S(y)$ then we may write
\begin{equation}\label{x2}
U\setminus f_S(ax + (1 - a)y) \subseteq U\setminus f_S(y).
\end{equation}
If $f_S(x)\subset f_S(y)$ then we may write
\begin{equation}\label{x3}
U\setminus f_S(ax + (1 - a)y) \subseteq U\setminus f_S(x).
\end{equation}
From (\ref{x2}) and (\ref{x2}), we have
\begin{equation}\label{x3}
U\setminus f_S(ax + (1 - a)y) \subseteq \{U\setminus f_S(x) \cup
U\setminus f_S(y)\}.
\end{equation}
So, $S^c$ is a concave soft set.
\begin{ttheorem}
$S^c$ is a convex soft set when S is a concave soft sets.
\end{ttheorem}
\emph{\textbf{Proof:}} Suppose that there exist $x, y\in E$, $a\in
I$ and S be a concave soft set.

Then, since S is concave,
\begin{equation}\label{y1}
f_S(ax + (1 - a)y) \subseteq f_S(x)\cup f_S(y)
\end{equation}
or
\begin{equation}\label{y2}
U\setminus f_S(ax + (1 - a)y) \supseteq U\setminus \{f_S(x)\cup
f_S(y)\}
\end{equation}
If $f_S(x)\supset f_S(y)$ then we may write
\begin{equation}\label{y3}
U\setminus f_S(ax + (1 - a)y) \supseteq U\setminus f_S(x).
\end{equation}
If $f_S(x)\subset f_S(y)$ then we may write
\begin{equation}\label{y4}
U\setminus f_S(ax + (1 - a)y) \supseteq U\setminus f_S(y).
\end{equation}
From (\ref{y3}) and (\ref{y4}), we have
\begin{equation}\label{x3}
U\setminus f_S(ax + (1 - a)y) \supseteq \{U\setminus f_S(x) \cap
U\setminus f_S(y)\}.
\end{equation}
So, $S^c$ is a convex soft set.
\begin{ttheorem} S is a concave soft set on E iff for every $\beta \in [0,1]$ and $\alpha \in P(U)$, $S^\alpha$ is a
concave set on E.
\end{ttheorem}
\emph{\textbf{Proof:}} $\Rightarrow$ Assume that S is a concave soft
set. If $x_1, x_2\in E$ and $\alpha \in P(U)$, then
$f_S(x_1)\supseteq \alpha$ and $f_S(x_2)\supseteq \alpha$. It
follows from the concavity of S that
$$
f_S(\beta x_1 +(1-\beta)x_2)\subseteq f_S(x_1)\cup f_S(x_2)
$$
and thus $S^\alpha$ is a concave set.

$\Leftarrow$ Assume that $S^\alpha$ is a concave set for every
$\beta\in [0, 1]$. Especially, for $x_1, x_2\in E,$ $S^\alpha$ is
concave for $\alpha =f_S(x_1)\cup f_S(x_2)$.

Since $f_S(x_1) \supseteq\alpha$ and $f_S(x_2) \supseteq\alpha$, we
have $x_1 \in S^\alpha$ and $x_2 \in S^\alpha$, whence $\beta x_1
+(1-\beta)x_2\in S^\alpha$. Therefore, $f_S(\beta x_1
+(1-\beta)x_2)\subseteq\alpha=  f_S(x_1)\cup f_S(x_2)$, which
indicates S is a concave soft set on X.
\section{Conclusion}
In the literature, convex fuzzy sets has been introduced widely by
many researchers. In this paper, we defined convex soft sets and
concave soft soft sets and give some properties. It can extend to
strictly convex and strongly convex soft set by using strictly
convex and strongly convex set which are useful in fuzzy
optimization.

It is worth pointing out that the results obtained here for convex
soft sets can be generalized to strictly convex and strongly convex
soft set by a similar way.

Also we will try to explore characterizations of convex soft sets to
optimization in the future.

The theory may be applied to many fields and more comprehensive in
the future to solve the related problems, such as; pattern
classification, operation research, decision making, optimization
problem, and so on.


\end{document}